\documentclass[11pt,reqno]{amsart}

\newtheorem{proposition}{Proposition}[section]

\newtheorem{corollary}[proposition]{Corollary}
\newtheorem{theorem}[proposition]{Theorem}

\theoremstyle{definition}
\newtheorem{definition}[proposition]{Definition}

\newtheorem{remark}[proposition]{Remark}

\newcommand{\thlabel}[1]{\label{th:#1}}

\newcommand{\selabel}[1]{\label{se:#1}}
\newcommand{\seref}[1]{Section~\ref{se:#1}}

\newcommand{\prlabel}[1]{\label{pr:#1}}
\newcommand{\prref}[1]{Proposition~\ref{pr:#1}}
\newcommand{\colabel}[1]{\label{co:#1}}
\newcommand{\coref}[1]{Corollary~\ref{co:#1}}
\newcommand{\relabel}[1]{\label{re:#1}}

\newcommand{\delabel}[1]{\label{de:#1}}

\newcommand{\eqlabel}[1]{\label{eq:#1}}
\newcommand{\equref}[1]{(\ref{eq:#1})}

\def\ZZ{{\mathbb Z}}

\newcommand{\Cc}{\mathcal{C}}

\def\*C{{}^*\hspace*{-1pt}{\Cc}}
\def\text#1{{\rm {\rm #1}}}

\def\catk{{{}_\Lambda\mathcal{C}}}
\def\Catk{{\mathcal{C}_\Lambda}}

\usepackage{amssymb}

\begin{document}

\title[Bicrossed products for finite groups] {Bicrossed products for finite groups}

\dedicatory{Dedicated to Freddy Van Oystaeyen  on the occasion of
his 60th birthday.}

\author{A. L. Agore}\thanks{The work of A. L. Agore, B. Ion and G. Militaru where
partially supported by CNCSIS grant 24/28.09.07 of PN II "Groups,
quantum groups, corings and representation theory". The work of B.
Ion was partially supported by NSF grant DMS-0536962. }
\author{A. Chirv\u asitu}
\address[A.L.A., A.C., G.M.]{Faculty of Mathematics and Computer Science, University of Bucharest, Str.
Academiei 14, RO-70109 Bucharest 1, Romania}
\email{ana.agore@fmi.unibuc.ro, chirvasitua@gmail.com,
gigel.militaru@fmi.unibuc.ro}
\author{B. Ion}
\address[B.I.]{Department of Mathematics, University of Pittsburgh,
Pittsburgh, PA 15260 }
\author{G. Militaru}


\begin{abstract}{ We investigate one question regarding bicrossed products of
finite groups which we believe has the potential of being
approachable for other classes of algebraic objects (algebras,
Hopf algebras). The problem is to classify the groups that can be
written as bicrossed products between groups of fixed isomorphism
types. The groups obtained as bicrossed products of two finite
cyclic groups, one being of prime order, are described.}
\end{abstract}

\maketitle

\section*{Introduction}

The bicrossed product construction is a generalization of the
semidirect product construction for the case when neither factor
is required to be normal: a group $E$ is the internal bicrossed
product of its subgroups $H$ and $G$ if $HG=E$ and their
intersection is trivial. Groups with this property (but allowing
for nontrivial intersection) have been in the literature for a
quite long time under the terminology {\sl permutable groups}
\cite{Maillet, Ore} or groups that admit an {\sl exact
factorization} (see e.g. \cite{Gi, WW}).

The bicrossed product construction itself is due to Zappa
\cite{zappa}. It was rediscovered by Sz\' ep \cite{Szep} and yet
again by Takeuchi \cite{Takeuchi}. The terminology {\sl bicrossed
product} is taken from Takeuchi,  other terms referring to this
construction used in the literature are {\sl knit product} and
{\sl Zappa-Sz\' ep product.} Bicrossed product constructions were
subsequently introduced and studied  for other structures:
algebras, Hopf algebras, Lie algebras, Lie groups, locally compact
quantum groups, groupoids. For Hopf algebras, in particular,
structural results are still missing and objects obtained from
such constructions form a considerable proportion of the known
examples (see e.g. \cite{CIMZ}). Assume for simplicity that  $k$
is a field of characteristic zero. Let $E$ be a finite group that
is a bicrossed product of the groups $H$ and $G$. A noncommutative
noncocommutative Hopf algebra $k[H]^* \# k[G]$ that is both
semisimple and cosemisimple can be constructed \cite{Takeuchi}.
This is the easiest way to construct semisimple cosemisimple
finite dimensional Hopf algebras. For this reason we decided to
investigate some aspects of the bicrossed product construction in
its original finite group setting.

Our main question, going back to Ore \cite{Ore} asks for the
description of all groups which arise as bicrossed products of two
fixed groups. Little progress has been made on this question. In
this respect we would like to mention the result of Wielandt
\cite{Wielandt} establishing that from two finite nilpotent groups
of coprime orders one always obtains a solvable group and the work
of Douglas \cite{Douglas} on finite groups expressible as
bicrossed products of two finite cyclic groups. Finding all
matched pairs between two finite cyclic groups seems to be still
an open question, even though J. Douglas \cite{Douglas} has
devoted four papers and over two dozen theorems to the subject. In
fact, solving this problem does not provide an answer to the
classification of all associated bicrossed products and does not
indicate whether a bicrossed product could not be obtained more
easily as a semidirect product. In \seref{3} we will give a
complete answer to this question for the case of two finite cyclic
groups, one of them being of prime order. As it turns out, if a
group is isomorphic to a bicrossed product of two finite cyclic
groups, one of them being of prime order then it is isomorphic to
a semidirect product between the same cyclic groups. We would like
also to mention some interesting recent investigations
\cite{Jara}, \cite{Pena} into the corresponding question at the
level of algebras.


\section{Prelimaries}

\subsection{Definitions and notation}\selabel{1}\selabel{1.2}
Let us fix the notation that will be used throughout the paper.
Let $H$ and $G$ be two groups and $\alpha : G \times H \rightarrow
H$ and $\beta : G \times H \rightarrow G$ two maps. We use the
notation
$$\alpha (g, h) = g\triangleright h \quad \text{and} \quad
\beta (g, h) = g\triangleleft h$$ for all $g\in G$ and $h\in H$.
The map $\alpha$ (resp. $\beta$) is called trivial if
$g\triangleright h = h$ (resp. $g\triangleleft h = g$) for all
$g\in G$ and $h\in H$. If $\alpha : G \times H \rightarrow H$ is
an action of $G$ on $H$ as group automorphisms we denote  by $H
\rtimes_{\alpha} G$ the semidirect product of $H$ and $G$: $H
\rtimes_{\alpha} G = H\times G$ as a set with the multiplication
given by
$$
(h_1,\, g_1)\cdot (h_2, \, g_2) : = \bigl( h_1 (g_1\triangleright
h_2), \, g_1 g_2 \bigl)
$$
for all $h_1$, $h_2 \in H$, $g_1$, $g_2 \in G$.

The opposite group structure on $H$ will be denoted by $H^{\rm
op}$ : $H^{\rm op} = H$ as a set with the multiplication $h_1
\cdot_{\rm op} h_2 = h_2 h_1$ for all $h_1$, $h_2 \in H$.

\begin{definition}\delabel{knitsyst}
A \textit{matched pair} of groups is a quadruple $\Lambda = (H, G,
\alpha, \beta)$ where $H$ and $G$ are groups, $\alpha : G \times H
\rightarrow H$ is a left action of the group $G$ on the set $H$,
$\beta : G \times H \rightarrow G$ is a right action of the group
$H$ on the set $G$ such that the following compatibility
conditions hold:
\begin{equation}\eqlabel{KS4}
g \triangleright (h_1 h_2) = (g\triangleright h_1)
((g\triangleleft {h_1} )\triangleright h_2)
\end{equation}
\begin{equation}\eqlabel{KS3}
(g_1 g_2)\triangleleft h = (g_{1}\triangleleft ({g_2
\triangleright h})) (g_2\triangleleft h)
\end{equation}
for all $g$, $g_1$, $g_2 \in G$, $h$, $h_1$, $h_2\in H$.

A morphism $\varphi:(H_1, G_1, \alpha_1, \beta_1)\to (H_2, G_2,
\alpha_2, \beta_2)$ between two matched pairs consists of a pair
of group morphisms  $\varphi_H: H_1\to H_2$, $\varphi_G:G_1\to
G_2$ such that
$$\varphi_H\circ\alpha_1=\alpha_2\circ (\varphi_G\times\varphi_H),
\quad \varphi_G\circ\beta_1=\beta_2\circ (\varphi_G\times\varphi_H)$$

\end{definition}

\begin{remark}
Let $\Lambda = (H, G, \alpha, \beta)$ be a matched pair of groups.
Then
\begin{equation}\eqlabel{KS5}
g\triangleright 1 = 1 \qquad {\rm and} \qquad 1\triangleleft h = 1
\end{equation}
for all $g\in G$ and $h\in H$.
\end{remark}

Let $H$ and $G$ be groups and $\alpha : G \times H \rightarrow H$
and $\beta : G \times H \rightarrow G$ two maps. Let $H\,
{}_{\alpha}\!\! \bowtie_{\beta} \, G = H\bowtie \, G : = H\times
G$ as a set with an binary operation defined by the formula:
\begin{equation}\eqlabel{knit4}
(h_1,\, g_1)\cdot (h_2, \, g_2) : = \bigl( h_1 (g_1\triangleright
h_2), \, (g_1\triangleleft {h_2}) g_2 \bigl)
\end{equation}
for all $h_1$, $h_2 \in H$, $g_1$, $g_2 \in G$.

The main motivation behind the definition of matched pair is the
following result (we refer to
 \cite{Takeuchi} or \cite[section IX.1]{Kassel} for the proof).

\begin{theorem}\thlabel{kn6}
Let $H$ and $G$ be groups and $\alpha$ and $\beta$ two maps as
above. Then $H\, {}_\alpha\!\! \bowtie_{\beta}\, G$  is a group
with unit $(1,1)$ if and only if $(H, G, \alpha, \beta)$ is a
matched pair. Moreover, a morphism between two matched pairs
induces a morphism between the corresponding groups.
\end{theorem}
If  $(H, G, \alpha, \beta)$ is a matched pair the group $H\bowtie
\, G$ is called the {\sl bicrossed product} (or the {\sl Zappa-Sz\'
ep product}) of $H$ and $G$. The inverse of an element of the
group $H\bowtie \, G$ is given by the formula
\begin{equation}\eqlabel{1.1.65k}
(h, g)^{-1} = \Bigl( g^{-1}\triangleright h^{-1}, \,  \bigl(
g\triangleleft {(g^{-1}\triangleright h^{-1})} \bigl)^{-1} \Bigl)
\end{equation}
for all $h\in H$ and $g\in G$. Also, remark that $H \times
\{1\}\cong H$ and $\{1\} \times G \cong G$ are subgroups of
$H\bowtie \, G$ and every element $(h, g)$ of $H\bowtie \, G$ can
be written uniquely as a product of an element of $H \times \{1\}$
and of an element of $\{1\} \times G$ as follows:
\begin{equation}\eqlabel{2.2.gen}
(h,g) = (h,1) \cdot (1,g)
\end{equation}
Conversely, one can see that this observation characterizes the
bicrossed product. Again, we refer to \cite{Kassel, Takeuchi} for the details.

\begin{theorem}\thlabel{kn7}
Let $E$ be a group $H$, $G\leq E$ be subgroups such that any
element of $E$  can be written uniquely as a product of an element
of $H$ and an element of $G$. Then there exists a matched pair
$(H, G, \alpha, \beta)$ such that
$$
\theta : H\bowtie \, G \ \rightarrow E, \qquad \theta (h, g) = hg
$$
 is group isomorphism.
\end{theorem}

The maps $\alpha$ and $\beta$ play in fact a symmetric role.
\begin{proposition}\prlabel{2.2.5n}
Let $\Lambda = (H, G, \alpha, \beta)$ be a matched pair of groups.
Then
\begin{enumerate}
\item[(i)] $\tilde{\Lambda} = (G, H, \tilde{\alpha}, \tilde{\beta}
)$, where $\tilde{\alpha}$ and $ \tilde{\beta}$ are given by
\begin{equation}\eqlabel{2.2.6n}
\tilde{\alpha} : H\times G \rightarrow G, \qquad \tilde{\alpha}
(h, g) = \Bigl(\beta (g^{-1}, h^{-1})\Bigl)^{-1}
\end{equation}
\begin{equation}\eqlabel{2.2.7n}
\tilde{\beta} : H\times G \rightarrow H, \qquad \tilde{\beta} (h,
g) = \Bigl(\alpha (g^{-1}, h^{-1})\Bigl)^{-1}
\end{equation}
for all $h\in H$ and $g\in G$ is a matched pair of groups.

\item[(ii)] The map
$$
\chi : \Bigl( H\, _{\alpha}\!\! \bowtie_{\beta} \, G   \Bigl)^{\rm op}
\rightarrow G\, _{\tilde\alpha}\!\! \bowtie_{\tilde\beta} \, H,
\qquad \chi (h,g) = (g^{-1}, h^{-1})
$$
is a group isomorphism. In particular,
$$
\xi :  H\, _{\alpha}\!\! \bowtie_{\beta} \, G  \rightarrow G\, _{\tilde\alpha}\!\! \bowtie_{\tilde\beta} \, H, \qquad \xi (h, g) =
\Bigl (g\triangleleft {(g^{-1}\triangleright h^{-1})},
\bigl(g^{-1}\triangleright h^{-1}\bigl)^{-1} \Bigl)
$$
is a group isomorphism.
\end{enumerate}
\end{proposition}

\begin{proof}  The proof is a straightforward verification.
\end{proof}

\begin{remark}\relabel{2.4.90}
Let $H$ and $G$ be two groups as above and let $\beta : G \times H
\rightarrow G$ be the trivial action. Then $(H, G, \alpha, \beta)$
is a matched pair if and only if the map $\alpha : G \times H
\rightarrow H$ is an action of $G$ on $H$ as group automorphisms.
In this case the bicrossed product is the semidirect product $H
\rtimes_{\alpha} G$.

Assume now that the map $\alpha$ is the trivial action. We obtain
from \equref{2.2.7n} that $\tilde{\beta}$ is trivial.  Keeping in
mind that the bicrossed product $G\, _{\tilde\alpha}\!\!
\bowtie_{\tilde\beta} \, H$ with the trivial ${\tilde{\beta}}$ is
a semidirect product we can invoke \prref{2.2.5n} (ii) to conclude
that
 $H\bowtie \, G \cong G\rtimes_{\tilde\alpha} H$.
\end{remark}

\subsection{Universality properties}\selabel{2-3}

Let $\Lambda = (H, G, \alpha, \beta)$ be a matched pair of groups.
We associate to $\Lambda$ two categories such that the bicrossed
product of $H$ and $G$ becomes an initial object in one of them
and a final object in the other.

Define  the category $\catk$ as follows: the objects of $\catk$ are pairs $(X, (u,v))$
where $X$ is a group, $u: H \rightarrow X$, $v: G \rightarrow X$
are group morphisms such that:
\begin{equation}\eqlabel{2.2.17}
v(g)u(h)= u(g \triangleright h)v (g \triangleleft h)
\end{equation}
for all $g\in G$, $h\in H$. A morphism  in $\catk$ $$f: (X, (u,v))
\rightarrow (X', (u',v'))$$ is a morphism of groups $f: X
\rightarrow X'$ such that $f \circ u = u'$ and $f\circ v = v'$. It
can be checked that $(H\bowtie \, G, (i_H, i_G))$ is an object in
$\catk$, where $i_H$ and $i_G$ are the canonical inclusions of $H$
and $G$ inside their bicrossed product.

Define the category $\Catk$ as follows: the objects of
$\Catk$ are pairs $(X, (u,v))$ where $X$ is a group, $u: X
\rightarrow H$, $v: X \rightarrow G$ are two maps such that the
following two compatibility condition holds:
\begin{equation}\eqlabel{2.2.18}
u(xy) = u (x) \bigl( v (x)\triangleright  u(y) \bigl), \qquad
v(xy) = \bigl( v(x) \triangleleft u(y)\bigl) v(y)
\end{equation}
for all $x$, $y \in X$. A morphism  in $\Catk$ $$f: (X, (u,v))
\rightarrow (X', (u',v'))$$ is a morphism of groups $f: X
\rightarrow X'$ such that $u' \circ f = u$ and $v'\circ f = v$. It
can be checked that  $\bigl(H\bowtie \, G,  \, (\pi_H, \pi_G)
\bigl)$ is an object in $\Catk$, where $p_H$ and $p_G$ are the
canonical projections from the bicrossed product to $H$ and $G$.

\begin{proposition}\prlabel{2.2.10}
Let $\Lambda = (H, G, \alpha, \beta)$ be a matched pair of groups.
Then
\begin{enumerate}
\item[(i)] $\bigl( H\bowtie \, G, (i_H, i_G) \bigl)$ is an initial
object of $\catk$.

\item[(ii)]  $\bigl(H\bowtie \, G,  \, (\pi_H, \pi_G) \bigl)$ is a
final object of $\Catk$.
\end{enumerate}
\end{proposition}

\begin{proof}
(i) Let $(X, (u,v))\in \catk$. We have to prove that there exists
a unique morphism of groups $w:  H\bowtie \, G \rightarrow X$ such
that $w \circ i_H = u$ and $w\circ i_G = v$.

Assume that $w$ satisfies this condition. Then using
\equref{2.2.gen} we have:
\begin{eqnarray*}
w ((h,g))  &=& w ((h, 1) \cdot (1,g))= w ((h,1))
w((1,g))\\
&=& (w \circ i_H) (h) (w \circ i_G)(g) = u (h)v(g)
\end{eqnarray*}
for all $h\in H$ and $g\in G$ and this proves that $w$ is unique.

If we define
$$
w : H \bowtie \, G \rightarrow X, \qquad w (h,g) = u(h)v(g)
$$
then
\begin{eqnarray*}
\hspace*{-2cm} w((h_1,g_1)\cdot (h_2,g_2)) &=& w (h_1
(g_1\triangleright
h_2), (g_1 \triangleleft h_2) g_2)\\
&=& u(h_1) u(g_1\triangleright  h_2) v(g_1\triangleleft  h_2 ) v (g_2)\\
&
{=}& u(h_1)v(g_1)u(h_2)v(g_2) \\ &=&
w((h_1, g_1)) w ((h_2, g_2))
\end{eqnarray*}
showing that $w$ is a morphism of groups.

Part (ii) follows by a similar argument.
\end{proof}

Straightforward from \prref{2.2.10} we obtain the description of
morphisms between a group and a bicrossed product.

\begin{corollary}\colabel{2.2.morf}
Let $E$ be a group and $(H, G, \alpha, \beta)$ a matched pair.
Then
\begin{enumerate}
\item[(i)] $w : H \bowtie \, G \rightarrow E$ is a group morphism if
and only if there exist $u: H \rightarrow E$ and $v: G \rightarrow
E$ group morphisms such that
$$
v(g)u(h)= u(g \triangleright h)v (g \triangleleft h), \quad {\rm
and} \quad w(h,g) = u(h)v(g)
$$
for all $h\in H$ and $g\in G$. \item[(ii)] $w : E \rightarrow H \bowtie
\, G$ is a morphism of groups if and only if there exist $u : E
\rightarrow H$ and $v : E \rightarrow G$ two maps such that
$$
u(xy) = u (x) \bigl( v (x)\triangleright  u(y) \bigl), \qquad
v(xy) = \bigl( v(x) \triangleleft u(y)\bigl) v(y)
$$
and $w (x) = (u(x), v(x))$ for all $x$, $y\in E$.
\end{enumerate}
\end{corollary}

\begin{remark}
\coref{2.2.morf} can be used to describe all morphisms or
isomorphisms between two matched pairs $H\, _{\alpha}\!\!
\bowtie_{\beta} \, G$ and $H\, _{\alpha'}\!\! \bowtie_{\beta'} \,
G$. However, the descriptions are rather technical and we will not
include them here.
\end{remark}


\section{Bicrossed products between finite cyclic groups}\selabel{3}

As mentioned in the Introduction the question of describing all
groups which arise as bicrossed products of two given groups was
asked by Ore. The first and, by our knowledge, the only systematic
study of this kind, for groups which arise as bicrossed products
of two finite cyclic groups, was employed by J. Douglas in 1951.
In his first paper on the subject \cite[pag. 604]{Douglas} Douglas
formulates the problem he wants to solve: describe all groups all
whose elements are expressible in the form $a^i b^j$ where $a$ and
$b$ are independent elements of order $n$ and, respectively, $m$.
What Douglas refers to as independent elements is in fact the
condition that the cyclic groups  generated by each of these
elements have trivial intersection. Therefore the problem can be
formulated as follows: describe all groups which arise as
bicrossed products of two finite cyclic groups.

In what follows $C_n$ and $C_m$ will be two cyclic groups of
orders $n$ and $m$. We denote by $a$ and $b$  a fixed generator of
$C_n$ and, respectively, $C_m$. For any positive integer $k$ we
denote by $\ZZ_k$ the ring of residue classes modulo $k$ and by
$S(\ZZ_k)$ the set of bijective functions from $\ZZ_k$ to itself.
Let $\alpha : C_m \times C_n \rightarrow C_n$, and $\beta: C_m
\times C_n \rightarrow C_m$ be two actions. They are completely
determined by two maps $\theta\in S(\ZZ_n)$ and $\phi\in \ZZ_m$
such that
$$
\alpha (b, a^x ) = a ^{\theta (x)}, \qquad \beta (b^y, a) =
b^{\phi(y)}
$$
for any $x \in \ZZ_n$, $y \in \ZZ_m$. Douglas \cite[Theorem
I]{Douglas} obtained necessary and sufficient conditions on the
pair of maps $(\theta,\phi)\in S(\ZZ_n)\times S(\ZZ_m)$ which
would make $(C_n, C_m, \alpha, \beta)$ a matched pair. The
functions $\theta$ and $\phi$ satisfying his conditions were
called {\sl conjugate special substitutions}.

Finding all pairs of conjugate special substitutions (or,
equivalently, all matched pairs between two finite cyclic groups)
seems to be still an open question. Furthermore, solving this
problem does not provide an answer to the classification problem
for the associated bicrossed products. In particular, it does not
indicate whether a bicrossed product could not be obtained more
easily as a semidirect product.

We investigate the structure of bicrossed products of two finite
cyclic groups, one of which has prime order. We shall prove now
our main result.

\begin{theorem}\thlabel{alex}
Let $p$ be a prime number, and $m\ge 2$ a positive integer. A
bicrossed product between two cyclic groups of orders $p$ and
respectively $m$ is isomorphic to a semidirect product between
cyclic groups of the same orders $p$ and $m$.
\end{theorem}

\begin{proof}
Let $G$ be a group, $C_p$ and $C_m$ two fixed cyclic subgroups of
$G$, of orders $p$ and $m$ respectively, such that $G = C_p\bowtie
C_m$. Let $a$, $b$ be generators for $C_p$ and $C_m$ respectively.
We can certainly assume that $C_p$ and $C_m$ are not normal in
$G$.

Let $H=C_m\cap aC_ma^{-1}$. A subgroup of a finite cyclic group is
uniquely determined by its order so $aHa^{-1}=H$. Since $C_m$ and
$aC_ma^{-1}$ are abelian, $H$ is central in both. Also note that
$C_m\ne aC_ma^{-1}$, because we assumed $C_m$ is not normal in
$G$. The centralizer $C(H)$ of $H$ in $G$ must then strictly
contain $C_m$. Since $C_m$ has index $p$ in $G$, it is maximal. It
follows that $C(H)=G$, that is $H$ is a central subgroup of $G$.

In fact, for any $1\leq k\leq p-1$, $ H$ is the intersection
between $C_m$ and $a^kC_ma^{-k}$. Indeed, $a^k$ generates $C_p$,
so any subgroup of $G$ normalized by $a^k$ must also be normalized
by $a$.

We are now going to work in the quotient group $G/H$. Let us
denote by $$\pi: G\to G/H$$ the canonical projection. For any
$1\leq k\leq p-1$, the groups $\pi(C_m)$ and $\pi(a^kC_ma^{-k})$
intersect trivially. We can now apply a theorem of Frobenius
\cite[Theorem 9.11 and Exercise 9.9]{rotman} for the group $G/H$
to conclude that the subgroup $\pi(C_p)$ is normal. This means
that $\pi({bab^{-1}})=\pi( a^t)$ for some $1\leq t\leq p-1$, or,
equivalently, that
$$bab^{-1}=a^tc$$ for some $c\in H$. Raising this identity to
power $p$, we find that $c^p=1$. We may assume that $c\ne 1$,
otherwise $C_p$ would be normal in $G$. Also, $t\ne 1$ otherwise
$a^{-1}ba=cb\in C_m$, which means that $C_m$ is normal in $G$,
contrary to our assumption. Hence, we can find an integer $u$ such
that $u(t-1)\equiv 1~(mod~p)$.

With the notation $\tilde a=ac^u$, we then have $b\tilde
ab^{-1}=\tilde a^t$. The proof is now finished: the subgroup of
$G$ generated by $\tilde a$ is normal in $G$, has order $p$ and
intersects $C_m$ trivially.
\end{proof}

\textbf{ACKNOWLEDGMENTS:} We thank the referee for detailed
suggestions that have helped us improve this paper.



\begin{thebibliography}{99}

\bibitem{CIMZ}
S. Caenepeel, B. Ion, G. Militaru, and  S. Zhu, The factorization
problem and the smash biproduct of algebras and coalgebras. {\sl
Algebr. Represent. Theory} {\bf 3} (2000),  no. 1, 19--42.

\bibitem{Douglas}
J. Douglas, On finite groups with two independent generators. I,
II, III, IV. {\sl Proc. Nat. Acad. Sci. U. S. A.} {\bf 37} (1951),
604--610, 677--691, 749--760, 808--813.

\bibitem{Gi} M. Giudici, Factorisations of sporadic simple groups.
{\sl J. Algebra}  {\bf 304}  (2006), no. 1, 311--323.

\bibitem{Jara} P. Jara, J. L\'{o}pez Pe\~ na, G. Navarro, D.
Stefan, On the classification of twisting maps between $K^n$ and
$K^m$, arXiv:0805.2874.

\bibitem{Kassel}  C. Kassel, Quantum groups. {\sl Graduate Texts in Mathematics} {\bf 155}.
Springer-Verlag, New York, 1995.

\bibitem{Pena} J. L\'{o}pez Pe\~ na and G. Navarro, On the classification and properties
of noncommutative duplicates. {\sl $K$--Theory} {\bf 38} (2008),
no. 2, 223--234.

\bibitem{Maillet}
E. Maillet, Sur les groupes \' echangeables et les groupes d\'
ecomposables. {\sl Bull. Soc. Math. France}  {\bf 28}  (1900),
7--16.

\bibitem{Ore} O. Ore, Structures and group theory. I. {\sl Duke Math. J.}
{\bf 3} (1937), no. 2, 149--174.

\bibitem{re}
L. R\'{e}dei, Zur Theorie der faktorisierbaren Gruppen, I. {\sl
Acta Math. Acad. Sci. Hungar.} {\bf 1} (1950), 74--98.

\bibitem{rotman} J. Rotman, An introduction to the theory of groups. Fourth edition.
{\sl Graduate Texts in Mathematics} {\bf 148}. Springer-Verlag, New York,
1995.

\bibitem{Szep}
J. Sz\' ep, \" Uber die als Produkt zweier Untergruppen
darstellbaren endlichen Gruppen. {\sl Comment. Math. Helv.} {\bf
22} (1949), 31--33.

\bibitem{Takeuchi}
M. Takeuchi, {Matched pairs of groups and bismash products of Hopf
algebras}.  {\sl Comm.  Algebra}  {\bf 9} (1981), no. 8, 841--882.

\bibitem{zappa}
G. Zappa, {Sulla costruzione dei gruppi prodotto di due dati
sottogruppi permutabili tra loro}. {\sl Atti Secondo Congresso Un.
Mat. Ital., Bologna, 1940},  119--125. Edizioni Cremonense, Rome,
1942.

\bibitem{WW}
J. Wiegold and A. G.  Williamson, The factorisation of the
alternating and symmetric groups. {\sl Math. Z.} {\bf 175} (1980),
no. 2, 171--179.

\bibitem{Wielandt} H. Wielandt,
\" Uber das Produkt paarweise vertauschbarer nilpotenter Gruppen.
{\sl Math. Z.} {\bf 55} (1951) 1--7.


\end{thebibliography}
\end{document}